\documentclass[11pt]{article}
\usepackage{amsmath,amsthm,amsfonts,amscd,bezier}
\usepackage{amssymb}
\usepackage[usenames]{color}
\newtheorem{theorem}{Theorem}%[section]
\newtheorem{lemma}[theorem]{Lemma}

\newtheorem{example}[theorem]{Example}
\newtheorem{definition}[theorem]{Definition}
\newtheorem{corollary}[theorem]{Corollary}
\newtheorem{question}[theorem]{Question}
\newtheorem{remark}[theorem]{Remark}

\newtheorem{proposition}[theorem]{Proposition}
\newcommand{\pf}{\noindent{\bf Proof:\ \ }}
\newcommand{\cqd}{{\hfill $\rule{2mm}{2mm}$}\vspace{1cm}}

\begin{document}

\pagestyle{plain} \pagenumbering{arabic}
\title{The Semiring of Values of an Algebroid Curve}
\author{Carvalho, E. and Hernandes, M.E.
\thanks{The first author
was partially supported by CAPES and the second one by CNPq.}\ \
\thanks{Corresponding author: Hernandes, M. E.; email: mehernandes@uem.br }}
\date{ \ }
\maketitle

\begin{center}


key words: Algebroid Curves, Semiring, Standard Basis.
\end{center}

\begin{center} Abstract
\end{center}

{\small We introduce the semiring of values $\Gamma$ with respect to
the tropical operations associated to an algebroid curve. As a set,
$\Gamma$ determines and is determined by the well known semigroup of
values $S$ and we prove that $\Gamma$ is always finitely generated
in contrast to $S$. In particular, for a plane curve, we present a
straightforward way to obtain $\Gamma$ in terms of the semiring of
each branch of the curve and the mutual intersection multiplicity of
its branches. In the analytical case, this allows us to connect
directly the results of Zariski and Waldi that characterize the
topological type of the curve.}

\section{Introduction}

Let $\mathbb{K}$ be an algebraically closed field. We
denote by $\mathbb{K}[[\underline{X}]]$ the ring
$\mathbb{K}[[X_1,\ldots ,X_n]]$ of formal power series in the indeterminates $X_1, \ldots , X_n$ with coefficients in the field $\mathbb{K}$.

In this paper, an {\it algebroid curve} in the $n$-dimensional space
$\mathbb{K}^n$ ($n>1$) is a proper radical ideal
$Q=\bigcap_{i=1}^{r}P_i \subset \mathbb{K}[[\underline{X}]]$ such
that $\mathcal{O}_i =\frac{ \mathbb{K}[[\underline{X}]]}{P_i}$ has
Krull dimension one for each isolated prime $P_i$ with $i\in
\{1,\ldots ,r\}$. Each $P_i$ is called a {\it branch} of the curve
$Q$. We will assume that the curve $Q$ is {\it non-degenerate}, that
is, $\text{dim}_{\mathbb{K}}\frac{\mathcal{M}}{\mathcal{M}^2}=n$,
where
$\mathcal{M}$ denotes the maximal ideal of the local, complete and reduced ring $\mathcal{O} = \frac{\mathbb{K}[[\underline{X}]]}{Q}$.

The integral closure $\overline{\mathcal{O}}_i$ of $\mathcal{O}_i$
in its quotient field is a discrete valuation domain isomorphic to
the ring $\mathbb{K}[[t_i]]$ and we have the inclusions (via
isomorphism)
$$ \ \ \ \ \ \ \ \ \ \ \ \ \ \ \ \ \ \mathcal{O} \subseteq \bigoplus_{i=1}^r\mathcal{O}_i
\subseteq \overline{\mathcal{O}} =
\bigoplus_{i=1}^r\overline{\mathcal{O}}_i = \bigoplus_{i =
1}^{r}\mathbb{K}[[t_i]], \ \ \ \ \ \ \ \ \ \ \ \ (1)$$\label{aneisLocais} where $\overline{\mathcal{O}}$ denotes
the integral closure of $\mathcal{O}$ in its total ring of
fractions.

If
$v_i:\overline{\mathcal{O}}_i\rightarrow\overline{\mathbb{N}}:=\mathbb{N}\cup\{\infty\}$
denotes the normalized discrete valuation of
$\overline{\mathcal{O}}_i$, where $v_i(0)=\infty$, for all
$i=1,...,r$, then the set $S_i:=\{v_i(g);\
g\in\mathcal{O}_i\setminus\{0\}\}\subseteq \mathbb{N}$ is
classically called the {\it semigroup of values} of $\mathcal{O}_i$.
Given a nonzero divisor $g\in\mathcal{O}$, we define
$v(g):=(v_1(g),\ldots ,v_r(g))\in\mathbb{N}^r$, where $v_i(g)$ means
the value of the homomorphic image of $g\in\mathcal{O}$ in
$\mathcal{O}_i$. In this way, we obtain the {\it semigroup of
values} of $\mathcal{O}$:
$$S:=\{v(g);\ g\ \mbox{is a nonzero divisor in}\  \mathcal{O}\}\subseteq \bigoplus_{i=1}^r S_i\subseteq \mathbb{N}^r.$$

In \cite{castellanos} it is described a method to obtain $S_i$ for a
branch of a space curve and in \cite{Abramo-Marcelo} it is presented
algorithms to compute $S_i$ and sets of values for any
$\mathcal{O}_i$-modules in $\overline{\mathcal{O}}_i$ as well.

Despite the fact semigroups of irreducible curves are finitely
generated, the same is not true for curves with several branches.
For instance, the semigroup of $Q=\langle XY\rangle=\langle
X\rangle\cap\langle Y\rangle$ is
$S=\{(0,0)\}\cup\{(1,1)+(\alpha_1,\alpha_2);\
(\alpha_1,\alpha_2)\in\mathbb{N}^2\}$ and it does not admit a finite
set of generators as an additive semigroup.

For an analytic plane curve, i.e., $n=2$, given by $Q=\langle
f\rangle$, where $f = \prod_{i=1}^{r}f_i \in \mathcal{M} \subset
\mathbb{C}\{X, Y\}$, Zariski (see \cite{equising} and
\cite{zariski}) shows that the topological type of $f^{-1}(0)$ is
completely characterized by the semigroup $S_i$ of each branch
$\langle f_i\rangle$, with $1 \leq i \leq r$, and the mutual
intersection multiplicity
$I(f_j,f_k)=\text{dim}_{\mathbb{C}}\frac{\mathbb{C}\{X,Y\}}{\langle
f_j,f_k\rangle}$, with $1\leq j<k\leq r$. By the other hand, Waldi
in \cite{waldi} obtains the topological characterization of the germ
$f^{-1}(0)$ by means of the semigroup $S$. In this way, a natural
question is:
\begin{question} How to obtain $S$ by means of $S_i$ and $I(f_j,f_k)$, for $1 \leq i \leq r$ and $1\leq j<k\leq r$?\end{question}

For two plane branches, Garcia and Bayer (see \cite{garcia} and
\cite{bayer}) answer this question using the notion of {\it
maximal points} of $S$. For a plane curve with several branches
$Q=\cap_{i=1}^{r}\langle f_i\rangle$, Delgado in \cite{delgado}
determines $S$ using the {\it relative maximal points} of $S$ and
the semigroups of $\cap_{i=1,i\neq j }^{r}\langle f_i\rangle$ for
all $1\leq j\leq r$.

%-novo
In \cite{AGMT} the authors consider (good) subsemigroups of $\mathbb{N}^r$ not necessarily associated to an algebroid curve, that is, under the arithmetical viewpoint. In that paper is described a finite subset $G$ of a good semigroup $S\subset \mathbb{N}^r$ (distinct of the maximal points consider by Delgado in \cite{delgado}) such that $G$ and the conductor of $S$ (see section 2) allow to determine $S$.

Given an algebroid curve $Q$ em $\mathbb{K}^n$, we propose to
consider the set
$$\Gamma=\{v(g);\ g\in\mathcal{O}\}\supset S,$$
where $v(0) = \underline{\infty}:= (\infty, \ldots , \infty)$.

Obviously, $(\Gamma,+)$ is a semigroup setting $\gamma +
\underline{\infty}=\underline{\infty}$ for all $\gamma\in\Gamma$. As
a semigroup, $\Gamma$ is the topological closure of $S$ in the
product topological space $\overline{\mathbb{N}}^r$, with
$\overline{\mathbb{N}}$ provided of the one point compactification
topology. This completion was considered for plane curves by Delgado
(see Section 2 of \cite{delgado}). In this paper, we will consider
the set $\Gamma$ equipped with the
tropical operations
$$\alpha \oplus\beta=\text{min}\{\alpha,\beta\}:=(\text{min}\{\alpha_1,\beta_1\},\ldots ,\text{min}\{\alpha_r,\beta_r\})\ \ \
\mbox{and} \ \ \  \alpha\odot\beta=\alpha+\beta,$$ for all $\alpha =
(\alpha_1, ..., \alpha_r)$ e $ \beta = (\beta_1, ..., \beta_r)$ em $
\Gamma$.

The main advantage of this approach is that $(\Gamma,\oplus,\odot)$
is a {\it finitely generated semiring} and, for plane curves, its
minimal set of generators allows us to connect $S$ with $S_i$ and $I(f_j,f_k)$ for
$1 \leq i \leq r$ and $1\leq j<k\leq r$, that is,
we obtain an answer to Question 1.

%-novo
Cotterill, Feital and Martins, in \cite{CFM}, used the tropical operations on the semigroup of values to obtain interesting results on singular rational curves in projective space.

\section{The Semiring of Values and Standard Bases}

Throughout this paper, we denote the set of indices $\{1,\ldots
,r\}$ by $I$.

Let $Q=\bigcap_{i\in I}P_i\subset \mathbb{K}[[\underline{X}]]$ be an
algebroid curve. As we remarked in the introduction, the integral
closure $\overline{\mathcal{O}}_i$ of the domain $\mathcal{O}_i
=\frac{\mathbb{K}[[\underline{X}]]}{P_i}$ is a discrete valuation
ring isomorphic to $\mathbb{K}[[t_i]]$, where $t_i$ is a uniformizing
parameter of $\overline{\mathcal{O}}_i$. In what follows, we
identify $X_j+P_i\in\mathcal{O}_i$ with its isomorphic image
$x_j(t_i)\in\mathbb{K}[[t_i]]$ and $\mathcal{O}_i$ with
$\mathbb{K}[[x_1(t_i),\ldots ,x_n(t_i)]]$. We call $(x_1(t_i),\ldots
,x_n(t_i))$ a parameterization of $P_i$.

If $v_i$ denotes the normalized valuation of
$\overline{\mathcal{O}}_i$, we have the additive semigroup
$$\Gamma_i= v_i(\mathcal{O}_i):=\{v_i(g_i)=ord_{t_i}(g_i);\
g_i\in\mathcal{O}_i\}\subseteq \overline{\mathbb{N}},$$ setting
$\gamma_i+\infty=\infty$ for all $\gamma_i\in\Gamma_i$ and $i\in I$.

Considering (\ref{aneisLocais}), we obtain the set of values of $\mathcal{O}$:
$$\Gamma=v(\mathcal{O}):=\{v(g):=(v_1(g),\dots ,v_r(g));\ g\in\mathcal{O}\}\subseteq\overline{\mathbb{N}}^r.$$

Notice that $S:=\Gamma\cap\mathbb{N}^r$ and $\Gamma$ determine each
other.

Consider a non empty subset $J=\{j_1,\ldots ,j_s\}$ of $I$. Setting
by $\Gamma_J$ the set of values of
$\mathcal{O}_J=\frac{\mathbb{K}[[\underline{X}]]}{\bigcap_{j\in
J}P_j}$, we denote by $Q^J$ the canonical image of the ideal
$\bigcap_{i\in I\setminus J}P_i+\bigcap_{j \in J}P_j\subset\mathbb{K}[[\underline{X}]]$
in $\mathcal{O}_J$ and by $v_J(Q^J)$ the $\Gamma_J$-monomodule
$\{(v_{j_1}(q),\dots ,v_{j_s}(q));\ q\in Q^J\}$. If $J=\{i\}$, we
put $Q^{\{i\}}=Q^{i}$ and $v_{\{i\}}(Q^{\{i\}})=v_i(Q^{i})$.

Since $\overline{\mathcal{O}}$ is an $\mathcal{O}$-module of finite
type, the conductor
$\mathcal{C}=(\mathcal{O}:\overline{\mathcal{O}})$ is an ideal of
$\overline{\mathcal{O}}$ and of $\mathcal{O}$ containing a nonzero
divisor and $\mathcal{C}=(t_1^{\sigma_1},\ldots
,t_r^{\sigma_r})\overline{\mathcal{O}}$. The element
$\sigma=(\sigma_1,\ldots ,\sigma_r)\in\Gamma$ is called the {\it
conductor} of $\Gamma$. As $\Gamma_i$ has a conductor, there exists
$\delta_i\in v_i(Q^{i})$ such that
$\delta_i+\overline{\mathbb{N}}\subseteq v_i(Q^{i})$ and $\delta_i$
is the smallest element in $v_i(Q^{i})$ with this property. D'Anna,
in \cite{danna} (Proposition 1.3), proves that $\sigma_i=\delta_i$
for all $i\in I$.

In \cite{Abramo-Marcelo} it is presented an algorithm to compute the
set of values for any finitely generated $\mathcal{O}_i$-module in
$\mathbb{K}[[t_i]]$. In particular, we can obtain $v_i(Q^{i})$ and
compute $\sigma_i$ for any $i\in I$.

\begin{remark} For a plane curve $Q=\langle \prod_{i\in I}f_i\rangle=\cap_{i\in I}\langle f_i\rangle$,
we have that $Q^{i}=\langle \prod_{j\in I\atop j\neq i}f_j\rangle$
and

\begin{center}

$v_i(Q^{i})=v_i\left( \prod_{j\in I\atop j\neq i}f_j\right)+\Gamma_i
= \sum_{j\in I\atop j\neq i}I(f_j,f_i)+\Gamma_i.$

\end{center}

\noindent In this way, we have the well know equality
$\sigma_i=\sum_{j\in I\atop j\neq i}I(f_j,f_i)+c_i$, where $c_i$ is
the conductor of $\Gamma_i$ that can be computed in terms of the
minimal set of generators of $\Gamma_i$.
\end{remark}

Let $\alpha = (\alpha_1, \ldots , \alpha_r)$ and $\beta = (\beta_1, \ldots , \beta_r)$ be elements in $\Gamma$. We have the
following properties:
\begin{itemize}
\item[a)] If $\alpha_i=0$ for some $i\in I$,
then $\alpha=\underline{0}:=(0,\ldots ,0)$.
\item[b)] If $\alpha_k=\beta_k < \infty$ for
some $k\in I$, then there exists $\gamma = (\gamma_1, \ldots , \gamma_r) \in\Gamma$ such that
$\gamma_i\geq\min\{\alpha_i,\beta_i\}$ for all $i\in I$ (the
equality holds if $\alpha_i\neq\beta_i$) and
$\gamma_k>\alpha_k=\beta_k$.
\item[c)] $\min\{\alpha,\beta\}:=(\min\{\alpha_{1},\beta_{1}\},\ldots
,\min\{\alpha_{r},\beta_{r}\})\in\Gamma$.
\end{itemize}

The last property allows us to consider $\Gamma$ equipped with the
tropical operations
$$\alpha \oplus\beta=\text{min}\{\alpha,\beta\}\ \ \ \ \ \text{and} \
\ \ \ \  \alpha\odot\beta=\alpha+\beta.$$ It is immediate that
$(\Gamma,\oplus,\odot)$ is a semiring.

\begin{definition} We call $(\Gamma,\oplus,\odot)$ the {\rm semiring of
values} associated to the curve $Q = \bigcap_{i\in I}P_i$.
\end{definition}

In the sequel we will show that $\Gamma$ is a finitely generated
semiring, that is, there exists a subset $\{\gamma_1,\ldots
,\gamma_m\}\subset\Gamma$ such that for any $\gamma\in\Gamma$ we can
write
$$\gamma=
\left (\gamma_1^{\alpha_{11}}\odot\ldots\odot \gamma_m^{\alpha_{1m}}\right )
\oplus\ldots \oplus \left (\gamma_1^{\alpha_{s1}}\odot\ldots \odot
\gamma_m^{\alpha_{sm}}\right )=\text{min}\left \{
\sum_{j=1}^{m}\alpha_{1j}\gamma_j,\ldots
,\sum_{j=1}^{m}\alpha_{sj}\gamma_j\right\}$$ with $\alpha_{ij}\in\mathbb{N},\  1\leq j\leq m$
and $1\leq i\leq s$, for some $s\leq r$ which depends on $\gamma$.

Remark that for $r=1$ we do not have novelty, so in what follows we
always consider $r\geq 2$.

Let $G = \{g_1, \ldots, g_m\}$ be a subset of
$\mathcal{M}\subset\mathcal{O}$. A {\it $G$-product} is an element
of the form
$$ G^{\alpha} = \prod_{j = 1}^{m}g_j^{\alpha_j},$$
with $\alpha = (\alpha_1, ..., \alpha_m) \in \mathbb{N}^m$.

From now on, given $\gamma \in \overline{\mathbb{N}}^r \setminus \{\underline{\infty}\}$ we put $I_{\gamma} = \{i \in I; \ \gamma_i \neq \infty \}$. If $g \in \mathcal{O}\setminus\{0\}$, then we denote $I_g = I_{v(g)}$.

\begin{definition}
Let $g$ be a nonzero element in $\mathcal{O}$ and consider $k\in
I_g$. An element $h \in \mathcal{O}$ is a $k$-reduction of $g$
modulo $G$ if there exist $c \in \mathbb{K}$ and a $G$-product
$G^{\alpha}$ such that
$$h = g - cG^{\alpha},$$
with $v_i(h) \geq v_i(g)$ for all $i \in I$ and $v_k(h)>v_k(g)$. We
say that $h$ is a reduction of $g$ modulo $G$ if $h$ is a
$k$-reduction of $g$ modulo $G$ for some $k \in I_g$.
\end{definition}

\begin{remark} Notice that if $g\in\bigcap_{i\in I\atop i\neq k}P_i\setminus P_k$
admits a reduction modulo $G$, then $g$ admits only a $k$-reduction
(because $I_g = \{k\}$) and there exists a $G$-product $G^{\alpha}$ such
that $v(g)=v(G^{\alpha})$. \label{obs}
\end{remark}

Now we are able to introduce the notion of a Standard Basis for
$\mathcal{O}$.

\begin{definition}
Let $G$ be a nonempty finite subset of $\mathcal{M}\setminus\{0\}$.
We say that $G$ is a {\rm Standard Basis} for $\mathcal{O}$ if every
nonzero element $g \in \mathcal{O}$ has a reduction modulo $G$.
\end{definition}

The next proposition allows us to present another characterization
of a Standard Basis for the local ring $\mathcal{O}$. For this
purpose we need the following lemma whose proof is analogous to
Lemma 1.8 of \cite{delgado_gorenstein}.

\begin{lemma} Let $J$ be a nonempty subset of $I$ and let
    $(\alpha_1,\ldots ,\alpha_r)\in \overline{\mathbb{N}}^r$ with
    $\sigma_j \leq \alpha_j < \infty$ for all $j\in J$. We have that
    $(\alpha_1,\ldots ,\alpha_r)\in \Gamma$ if and only if there exists
    $(\beta_1, ..., \beta_r)\in\Gamma$, where $\beta_i=\alpha_i$ for all
    $i\in I\setminus J$ and $\beta_j=\infty$ for all $j\in J$.
    \label{lema}
\end{lemma}

\begin{proposition}
Let $G$ be a nonempty and finite subset of $\mathcal{M}\setminus\{0\}$. The
following statements are equivalent:

\noindent \emph{(a)} Every nonzero element $g \in \mathcal{O}$ has a
$k$-reduction modulo $G$ for all $k \in I_g$.

\noindent \emph{(b)} Every nonzero element $g \in \mathcal{O}$ has a
reduction modulo $G$.
\end{proposition}\label{proposicaoReducao}

\pf It is sufficient to prove (b) $\Rightarrow$ (a). Suppose,
contrary to our claim, that there exist a nonzero element $g
\in\mathcal{O}$ and $k\in I_g$ such that $g$ admits a reduction modulo $G$,
but it does not have a $k$-reduction. We may assume that $g$ just has a
$j$-reduction for all $j$ in a nonempty subset $J$ of
$I_g\setminus\{k\}$.

Now consider an element $h \in \mathcal{O}$ obtained from $g$ via a
finite chain of reductions modulo $G$ such that $h$ does not have a
$j$-reduction modulo $G$ or $\sigma_j \leq v_j(h) < \infty$ for all $j \in
J$, where $\sigma = (\sigma_1, ..., \sigma_r)$ is the conductor of $\Gamma$. Notice that $h\neq 0$ and $v_i(h)=v_i(g)$ for all $i\in I\setminus
J$. In this way $h$ does not admit an $i$-reduction for any $i\in
I_g\setminus J$, otherwise the same would be true for $g$.

Let us consider $L\subseteq J$ such that $h$ does not admit an
$i$-reduction for all $i\in I_g\setminus L\supseteq I_g\setminus J$,
that is, $\sigma_l \leq v_l(h) < \infty$ for all $l\in L$. If $L=\emptyset$,
then $h$ does not have a reduction modulo $G$, which is a
contradiction. By the other hand, if $L\neq\emptyset$, Lemma \ref{lema} implies that there exists $h' \in
\mathcal{O}\setminus\{0\}$ such that $v_i(h') =v_i(h)$ for all $i\in
I\setminus L$ and $v_l(h')=\infty$ for all $l\in L$. But in this
way, $h'$ does not admit a reduction modulo $G$ and we obtain
a contradiction again. \cqd

As an immediate consequence of the concept of reduction and the
above proposition, we have the following characterization for a Standard Basis for $\mathcal{O}$.

\begin{corollary}
Let $G$ be a nonempty and finite subset of
$\mathcal{M}\setminus\{0\}$. The following statements are
equivalent:

\noindent \emph{(a)} $G$ is a Standard Basis for $\mathcal{O}$.

\noindent \emph{(b)} For every nonzero element $g \in \mathcal{O}$
and for some $k\in I_g$, there exists a $G$-product $G^{\alpha}$
(which depends on $k$) such that $v_i(g) \leq v_i(G^{\alpha})$ for all
$i \in I$ and $v_k(g)=v_k(G^{\alpha})$.

\noindent \emph{(c)} For every nonzero element $g \in \mathcal{O}$
and for all $k\in I_g$, there exists a $G$-product $G^{\alpha}$
(which depends on $k$) such that $v_i(g) \leq v_i(G^{\alpha})$ for all
$i \in I$ and $v_k(g)=v_k(G^{\alpha})$. \label{corolarioReducaoSB}
\end{corollary}

In \cite{Abramo-Marcelo}, the notion of Standard Basis was introduced
for branches and, in that case, its existence is
immediate. The next theorem guarantees the existence of a
Standard Basis for the local ring of any algebroid curve with several branches.

\begin{theorem}\label{base}
The local ring $\mathcal{O}$ admits a Standard Basis.
\end{theorem}
\pf Let $H$ be a subset of $\mathcal{O}$ satisfying $v(H)
=v(\mathcal{M})= \Gamma\setminus\{\underline{0}\}$ such that $v(h) \not\in
v(H\backslash \{h\})$ for all $h \in H$ and set $B_0:=\{h \in H;\
v_i(h) <\sigma_i\ \mbox{if}\ i\in I_h\}$.

For all $i\in I$, consider $B'_i, B''_i \subset \mathcal{O}$ such that the
homomorphic image of $B'_i$ and $B''_i$ in $\mathcal{O}_i$ are Standard Bases for $\mathcal{O}_i$ and $Q^i$ respectively, which can be computed as described in \cite{Abramo-Marcelo}. As the
homomorphic image of any finite subset $A$ of $\mathcal{O}$ such
that $v_i(A)=v_i(B''_i)$ is a Standard Basis for $Q^i$, we can take
$B''_i$ as a subset of the homomorphic image of $\bigcap_{j\in I\atop
j\neq i}P_j$ in $\mathcal{O}$, that is, $v_j(h)=\infty$ for all $h \in B''_i$ and
$j\in I\setminus\{i\}$. Setting $B_i = B'_i \cup B''_i$, we will show that the finite set $G=\bigcup_{i=0}^{r}B_i$ is a Standard
Basis for $\mathcal{O}$.

Let $g$ be a nonzero element in $\mathcal{O}$.

If $v_i(g)<\sigma_i$ for all $i \in I_g$, then there exists a $G$-product $G^{\alpha}$
(more specifically a $B_0$-product) such that $v(g)= v(G^{\alpha})$.

If $\sigma_k\leq v_k(g)$ for some $k \in I_g$, then $v_k(g)\in v_k(Q^k)$. As the
homomorphic image of $B'_k, B''_k\subset\mathcal{O}$ are Standard Bases for
$\mathcal{O}_k$ and $Q^k$ respectively, there exists a $G$-product $G^{\alpha}$ (indeed, a $B'_k$-product $B^{\beta}$ and $h_k \in B''_k$ with $G^{\alpha} = B^{\beta}h_k$) such that
$v_k(g)=v_k(G^{\alpha})$ and $v_i(g)\leq v_i(G^{\alpha})=\infty$ for
all $i\in I\setminus\{k\}$.

By the above corollary, we conclude that $G$ is a Standard Basis for $\mathcal{O}$.
\cqd

The above theorem allows us to conclude that the semiring of
values associated to an algebroid curve is finitely generated.

\begin{theorem}
The semiring $\Gamma$ is generated by $v(G)$, where $G$ is a Standard Basis for $\mathcal{O}$.
\end{theorem}
\pf Let $G = \{g_1,\ldots ,g_m\}$ be a Standard Basis for $\mathcal{O}$ with $v(g_j)=\gamma_j=(\gamma_{j1},\ldots ,\gamma_{jr})\in\overline{\mathbb{N}}^r$, for $1 \leq j \leq m$.

Initially notice that $\underline{0}=\sum_{j=1}^{m}0\cdot\gamma_j=
\gamma_1^{0}\odot ... \odot \gamma_m^{0}$.

By Remark \ref{obs}, for each $0 \neq h_k\in Q^k$, there exists a $G$-product $G^{\beta_k}$ such that
$v(h_k)=v(G^{\beta_k})$. In this way,
$$\underline{\infty}=v(h_1)+v(h_2)=v(G^{\beta_1})+v(G^{\beta_2})=v(G^{\alpha})=\sum_{j=1}^{m}\alpha_j\cdot\gamma_j = \gamma_1^{\alpha_1}\odot ... \odot\gamma_m^{\alpha_m},$$
where $\beta_1+\beta_2 = \alpha = (\alpha_1, ..., \alpha_m)$.

Now, given $\rho = (\rho_1, \ldots , \rho_r) \in \Gamma\setminus
\{\underline{0},\underline{\infty}\}$, there exists $g \in
\mathcal{M}\setminus\{0\}$ such that $\rho = v(g)$.

If $I_g=\{i_1,\ldots ,i_s\}$ then, by Corollary \ref{corolarioReducaoSB}, there exist $\alpha_{kj}\in\mathbb{N}$ with $1\leq k\leq s$ and $1 \leq j \leq m$ such that $$\rho_{i_k} = v_{i_k}(g_1^{\alpha_{k1}}\cdot \ldots \cdot g_m^{\alpha_{km}}) \ \text{and} \ \rho_{i} \leq v_{i}(g_1^{\alpha_{k1}}\cdot \ldots \cdot g_m^{\alpha_{km}}), \ \text{for all} \ i \in I\setminus\{i_k\}.$$

In this way, for $i_k\in I_g$, we have
$$\rho_{i_{k}} = \min\left\{\sum_{j=1}^{m}\alpha_{1j}\gamma_{ji_{k}}, \ldots , \sum_{j=1}^{m}\alpha_{sj}\gamma_{ji_{k}}\right\}.$$

Therefore,
$$ \rho = \min
\left\{\sum_{j=1}^{m}\alpha_{1j}\gamma_j, ...,
\sum_{j=1}^{m}\alpha_{sj}\gamma_j \right\} = (\gamma_1^{\alpha_{11}}\odot
... \odot \gamma_m^{\alpha_{1m}}) \oplus ... \oplus
(\gamma_1^{\alpha_{s1}}\odot ... \odot \gamma_m^{\alpha_{sm}}),$$ that is,
the semiring $\Gamma$ is finitely generated by $v(G)$. \cqd

\begin{remark} Consider the $\mathcal{M}$-adic topology on
$\mathcal{O}$ and let $G=\{g_1,\ldots ,g_m\}$ be a Standard Basis
for $\mathcal{O}$. Given $g\in\mathcal{O}\setminus\{0\}$, we have a
chain (possibly infinite) of reductions modulo $G$
$$h_0=g,\ \ \ \ h_i=h_{i-1}-c_iG^{\alpha_i},\ i > 0,$$
where $c_i\in\mathbb{K}$ and $G^{\alpha_i}$ is a $G$-product.

In case of an infinite chain of reductions, we get a sequence
$s_k=\sum_{i=1}^{k}c_iG^{\alpha_i},\ k\geq 1$ in $\mathcal{O}$. As we have $v(G^{\alpha_i})\neq v(G^{\alpha_j})$ for $i\neq j$, the set
$\{c_iG^{\alpha_i};\ i\geq 1\}$ is summable and the sequence $s_k$
is convergent in $\mathcal{O}$.

Since $G$ is a Standard Basis for $\mathcal{O}$, every
element $g\in\mathcal{O}\setminus\{0\}$ admits a chain of reductions
modulo $G$ to $0$, that is,
$g=\lim_{k\rightarrow\infty}\sum_{i=1}^{k}c_iG^{\alpha_i}$ or,
equivalently, $\mathcal{O}=\mathbb{K}[[g_1,\ldots ,g_m]]$.
\end{remark}

\section{Minimal Standard Bases and Irreducible Absolute Points}

The semigroup $S$ of an irreducible algebroid curve admits a minimal
set of generators $V$ in the sense that every system of generators
of $S$ contains $V$. In this way, it is natural analyze this
property for the semiring of any algebroid curve with
several branches that, in turn, is closely related with properties of a Standard Basis for
$\mathcal{O}$.

It is obvious that for any Standard Basis $G$ for $\mathcal{O}$ and for every nonzero element $g\in\mathcal{M}$, the set $G\cup\{g\}$ is also a Standard Basis for $\mathcal{O}$. So it is convenient to introduce the following definition.

\begin{definition}
Let $G$ be a Standard Basis for $\mathcal{O}$. We say that $G$ is
{\rm minimal} if for every $g \in G$ there does not exist a reduction of
$g$ modulo $G\setminus \{g\}$.
\end{definition}

In the next proposition we will prove that from a Standard Basis $G$
we can always get a minimal Standard Basis discarding elements $g\in
G$ that admit some reduction modulo $G\setminus\{g\}$. This will
guarantee the existence of a minimal Standard Basis for
$\mathcal{O}$.

\begin{proposition}
Let $G$ be a Standard Basis for $\mathcal{O}$. If $g\in G$ admits some reduction modulo $H = G\setminus\{g\}$, then $H$ is a Standard Basis for $\mathcal{O}$.
\end{proposition}
\pf Suppose that $g$ admits a $k$-reduction modulo $H$ for some $k
\in I_g$, that is, there exist $c_1 \in \mathbb{K}$ and an
$H$-product $H^{\alpha_1}$ such that $h = g - c_1H^{\alpha_1}$
satisfies $v_i(h) \geq v_i(g)$ for all $i \in I$ and $v_k(h) >
v_k(g)$. If $v_i(g) = v_i(H^{\alpha_1})$ for all $i \in I$, by
Corollary 9, we have that $g$ admits an $i$-reduction modulo $H$ for
all $i \in I_g$. Therefore, $H$ is a Standard Basis for
$\mathcal{O}$.

By the other hand, if there exists $j \in I_g$ such that $v_j(g) <
v_j(H^{\alpha_1})$, then $j\in I_h$ and $h$ admits a $j$-reduction
modulo $G$, since $G$ is a Standard Basis for $\mathcal{O}$.
Consequently, there exist $c_2 \in \mathbb{K}$ and a $G$-product
$G^{\beta}$ such that the element $h' = h - c_2G^{\beta} = g -
c_1H^{\alpha_1} - c_2G^{\beta}$ satisfies $v_i(h') \geq v_i(h) \geq
v_i(g)$ for all $i \in I$ and $v_j(h') > v_j(h) = v_j(g).$

In this way, we must have $G^{\beta}=g$ or $G^{\beta}$ is an
$H$-product $H^{\alpha_2}$.

If $G^{\beta} = g$, then $c_2 = 1$, $h' = -c_1H^{\alpha_1}$ e
$v_k(h') = v_k(H^{\alpha_1}) < v_k(h)$, which is a contradiction. It
follows that $G^{\beta}=H^{\alpha_2}$ and we obtain $v_i(g) \leq
v_i(h) \leq v_i(H^{\alpha_2})$ for all $i \in I$ and $v_j(g) =
v_j(H^{\alpha_2})$. By Corollary 9, we conclude that $g$ admits a
$j$-reduction modulo $H$. Hence, $H$ is a Standard Basis for
$\mathcal{O}$. \cqd

It is easy to see that the elements in a minimal Standard Basis have pairwise distinct values. Moreover, we have the following result.

\begin{proposition} If $G$ and $H$ are Standard Bases for $\mathcal{O}$ with
 $H$ minimal, then $v(H)\subseteq v(G)$. In particular, all the minimal Standard Bases for $\mathcal{O}$ have the same set of values.
\end{proposition}
\pf Let $H=\{h_1,\ldots ,h_s\}$ and $G=\{g_1,\ldots ,g_m\}$ be Standard Bases for $\mathcal{O}$ such that $H$ is minimal.

We will show that $v(h_l)\in v(G)$ for all $h_l\in H$.

Without loss of generality, we can consider $l = 1$. As $G$ is a
Standard Basis for $\mathcal{O}$, given $k\in I_{h_1}$ there exists a $G$-product
$g_1^{\alpha_1}\cdot\ldots\cdot g_m^{\alpha_m}$ (with $\alpha_j=0$ if
$v_k(g_j)=\infty$) such that
$$v_i(h_1)\leq v_i(g_1^{\alpha_1}\cdot\ldots\cdot g_m^{\alpha_m})\ \mbox{for all}\ i\in I \ \mbox{and} \ v_k(h_1)=v_k(g_1^{\alpha_1}\cdot\ldots\cdot g_m^{\alpha_m}).$$

By the other hand, for each $j\in\{1,\ldots ,m\}$, with $v_k(g_j)\neq\infty$, there exist an $H$-product $h_1^{\beta_{j1}}\cdot\ldots\cdot h_s^{\beta_{js}}$ such that
$$v_i(g_j)\leq v_i(h_1^{\beta_{j1}}\cdot\ldots\cdot h_s^{\beta_{js}})\ \mbox{for all}\ i\in I \ \mbox{and} \ v_k(g_j)=v_k(h_1^{\beta_{j1}}\cdot\ldots\cdot h_s^{\beta_{js}}).$$

But, in this way, we have
$$v_i(h_1)\leq v_i(h_1^{\sum_{j=1}^{m}\alpha_j\beta_{j1}}\cdot\ldots\cdot h_s^{\sum_{j=1}^{m}\alpha_j\beta_{js}})\ \mbox{for all}\ i\in I \ \mbox{and} \ v_k(h_1)=v_k(h_1^{\sum_{j=1}^{m}\alpha_j\beta_{j1}}\cdot\ldots\cdot h_s^{\sum_{j=1}^{m}\alpha_j\beta_{js}}).$$

In particular, we must have $\sum_{j=1}^{m}\alpha_j\beta_{j1}\leq 1$.

If $\sum_{j=1}^{m}\alpha_j\beta_{j1}=0$, then $h_1$ admits a $k$-reduction
modulo $H\setminus\{h_1\}$, which contradicts the fact that $H$ is a
minimal Standard Basis for $\mathcal{O}$.

It follows that $\sum_{j=1}^{m}\alpha_j\beta_{j1}=1$ and
$\sum_{j=1}^{m}\alpha_j\beta_{j2}=\ldots =\sum_{j=1}^{m}\alpha_j\beta_{js}=0$. So,
there exists $j_0\in\{1,\ldots ,m\}$ such that $\alpha_{j_0}=\beta_{j_01}=1,
\alpha_{j_0}\beta_{j_02}=\ldots =\alpha_{j_0}\beta_{j_0s}=0$ and, consequently,
$\beta_{j_0l}=0$ for $l = 2,\ldots ,s$. Then we obtain
$v_k(g_{j_0})=v_k(h_1)$ and $\alpha_j=0$ for all $j \neq j_0$. In
addition, $v_i(h_1)\leq v_i(g_{j_0})\leq v_i(h_1)$.

Therefore, $v(h_1)=v(g_{j_0})\in v(G)$. \cqd

By the above proposition, if $G$ is a minimal Standard Basis for
$\mathcal{O}$, then $v(G)$ is the unique minimal system of generators for
the semiring of values $\Gamma$.

In what follows we will continue to explore the relationship between
a Standard Basis $G$ for $\mathcal{O}$ and the semiring of values $\Gamma$.

\begin{definition} An element
$\gamma\in\Gamma\setminus\{\underline{0}\}$ is called {\rm
irreducible} if
$$\gamma=\alpha+\beta;\ \alpha,\beta\in\Gamma
\ \ \Rightarrow\ \ \alpha=\gamma\ \ \mbox{or}\ \ \beta=\gamma.$$
\end{definition}

Notice that the value of any element in a minimal Standard Basis is an irreducible element of the semiring. The algebraic
counterpart of this property is true as well, that is, every element
in a minimal Standard Basis $G$ is irreducible in $\mathcal{O}$.
Indeed, if $g = g_1g_2 \in G$, where $g_1,g_2\in \mathcal{M}$, then $\gamma
=v(g)=v(g_1)+v(g_2)=\alpha +\beta$, with $\alpha\neq\gamma\neq\beta$.

Given $\gamma\in\Gamma$ and a proper subset $J$ of $I_{\gamma}$, we
set:
$$F_{J}(\gamma)=\{\alpha\in\Gamma;\ \alpha_i>\gamma_i\ \mbox{for}\ i\in I_{\gamma}\setminus J\ \mbox{and}\ \alpha_j=\gamma_j\
\mbox{for}\ j\not\in I_{\gamma}\setminus J\}.$$

\begin{remark}\label{absoluto-reducao} If $v(g)=\gamma\neq\underline{\infty}$ and
$F_J(\gamma)\neq\emptyset$ for some proper subset $J$ of
$I_{\gamma}$, then for every $j \in J$ there exists a $j$-reduction of $g$ modulo a Standard Basis $G$. By the
other hand, if $F_J(\gamma)=\emptyset$ for all $\emptyset\neq
J\subset I_{\gamma}$, then the only possibility of reduction of $g$ modulo $G$
is $h=g-cG^{\alpha}$, with $v_i(h)>v_i(g)$ for all $i\in I_{\gamma}$,
i.e., $v(g)=v(G^{\alpha})$ .
\end{remark}

\begin{definition} We say that $\gamma\in\Gamma$ is an {\rm absolute (maximal) point} of $\Gamma$ if $F_{J}(\gamma)=\emptyset$ for every proper subset $J$ of
$I_{\gamma}$.
\end{definition}

Notice that if $I_{\gamma}$ has only one element, then there does not
exist a proper subset $J$ of $I_{\gamma}$ such that
$F_J(\gamma)\neq\emptyset$. In this way, vacuously, $\gamma$ is
considered an absolute point of $\Gamma$.

For $r=1$, the previous definition is equivalent to say that every element in $\Gamma$ is an absolute point and the minimal
system of generators of $\Gamma$ is precisely its set of irreducible absolute
points.

Now we will show that, similarily to the irreducible case, $\Gamma$ is a
semiring minimally generated by its irreducible absolute
points.

\begin{theorem}
Let $G$ be a Standard Basis for $\mathcal{O}$ such that its elements have pairwise distinct values. Then $G$ is minimal if and only if $v(G)$ is the set
of irreducible absolute points of $\Gamma$.
\end{theorem}
\pf Assume that $G=\{g_1,\ldots ,g_m\}$ is a minimal Standard Basis for $\mathcal{O}$ and consider
$g\in G$. In particular, $v(g)\neq\underline{0}$ is
irreducible.

If $v(g)$ is not an absolute point of $\Gamma$, then
there exist $h \in \mathcal{O}$ and a nontrivial subset $J$ of $I_g$
such that $v_i(g) < v_i(h)$ for all $i \in I_g \setminus J$ and
$v_j(g) = v_j(h)$ for all $j\not\in I_g\setminus J$ or, equivalently,
for all $j \in J \cup (I\setminus I_g)$. As $G$ is a Standard Basis, for
each $k \in J$ there exists a $G$-product $G^{\alpha}$ (depending on
$k$) such that $v_i(h) \leq v_i(G^{\alpha})$ for all $i \in I$ and
$v_k(h)=v_k(G^{\alpha})$. Hence, $v_i(g) \leq v_i(G^{\alpha})$ for
every $i \in I$ with $v_k(g) = v_k(G^{\alpha})$ and $v_j(g)<
v_j(G^{\alpha})$ for $j\in I_g\setminus J$. But, in this way,
$G^{\alpha}$ is a $G\setminus \{g\}$-product, that is, there exists
a reduction of $g$ modulo $G\setminus \{g\}$, an absurd because $G$
is minimal. Therefore, $v(g)$ is an irreducible absolute
point of $\Gamma$.

Now, let $\gamma=v(g)$ be an irreducible absolute point of
$\Gamma$. By Remark \ref{absoluto-reducao}, there exists a
$G$-product such that $v(g) = v(G^{\alpha})$. Furthermore, since the element
$\gamma$ is irreducible, we must have $G^{\alpha} = g_1^0 \cdot ...
\cdot g_j^1 \cdot ... \cdot g_{m}^0$, for some $1\leq j\leq m$.
Therefore, $\gamma = v(g_j) \in v(G)$.

Conversely, assume that $v(G)$ is the set of all irreducible absolute points of $\Gamma$ and let $g$ be an element of $G$. Remark \ref{absoluto-reducao} implies that $v(g) = v(G^{\alpha})$ for some $G$-product $G^{\alpha}$. Now, as the elements in $G$ have pairwise distinct values and $v(g)$ is irreducible we must have $G^{\alpha} = g$, that is,
$g$ does not have a reduction modulo $G\setminus \{g\}$. Hence,
$G$ is a minimal Standard Basis for $\mathcal{O}$.
\cqd

As an immediate consequence we have the following result.

\begin{corollary}
If $G$ is a finite subset of $\mathcal{M}$ such that $v(G)$ is precisely the set of irreducible absolute points of $\Gamma$, then $G$ is a minimal Standard Basis for $\mathcal{O}$.
\end{corollary}

We will present in \cite{carvalho} an algorithm that allows us to
obtain a Standard Basis for $\mathcal{O}$ and consequently the minimal
system of generators for the semiring of values $\Gamma$ of any algebroid curve in $\mathbb{K}^n$. However, for plane curves we can give a direct and more precisely description.

Let $Q = \cap_{i=1}^{r}\langle f_i\rangle$ be a plane
curve and let $S$ be its semigroup of values. Given $J\subseteq I$, if $\pi_J$ denotes the natural
projection of $\overline{\mathbb{N}}^r$ to the set of indices $J$,
then $\pi_{\{i\}}(S)=S_i$ and the number of absolute points of
$\pi_{\{j,k\}}(S)$ is precisely the intersection multiplicity
$I(f_j,f_k)$ (see \cite{garcia}). In this way, $S$ determines $S_i$ and $I(f_j,f_k)$ for $i\in I$ and $1\leq j<k\leq r$.

By the other hand, as we remarked in the introduction, Delgado in \cite{delgado}
characterizes $S$ in terms of the semigroup of
the curves $Q^j=\langle \prod_{i\in I\atop i\neq j}f_i\rangle$ and a
finite subset $R$ (the set of relative maximal points) of $S$ (see the
Generation Theorem in \cite{delgado}). In order to obtain $R$, it
is computed the set of irreducible absolute points of $S$, which
corresponds to the set $A$ of values of curves with maximal contact with
some branch of $Q$ and, using a symmetry property with respect to the
conductor of $S$, he gets a set that
contains $R$ and this set allows to apply the Generation Theorem. We notice that $A$ is precisely the set of irreducible absolute points of
$S=\Gamma\cap\mathbb{N}^r$ and it can be obtained by $S_i$ and
$I(f_j,f_k)$, for $ i\in I$ and $ 1\leq j<k\leq r$.

The next proposition describes the irreducible absolute points of $\Gamma$ with some coordinate equal to
$\infty$ in terms of $I(f_j,f_k)$ and, as the semiring $\Gamma$ and
the semigroup $S$ determine each other, provides an answer to Question 1.

\begin{proposition}
The set of irreducible absolute points of $\Gamma$ with some
infinite coordinate for a plane curve $Q = \cap_{i\in
I}\langle f_i\rangle$ is $\{v(f_i);\ i\in I\}$.
\end{proposition}
\pf If $\gamma =v(f_i)$ is not an absolute point of $\Gamma$ for
some $i\in I$, then there exists a nontrivial subset $J$ of
$I_{f_i}=I\setminus\{i\}$ such that $F_J(\gamma)\neq\emptyset$,
i.e., there exists $\beta\in\Gamma$ such that $\beta_j =\gamma_j$
for all $j\in J$, $\beta_i =\gamma_i=\infty$ and $\beta_k >\gamma_k$
for all $k\not\in J$ and $k\neq i$. As $\beta_i=\gamma_i=\infty$, we
conclude that $\beta=v(hf_i)$ for some $h\in\mathcal{O}$. By the other hand, the equality $\beta_j=\gamma_j$, for all $j\in J$, implies $v(h)=\underline{0}$, but, in this way, $\beta =\gamma$, which is an absurd
because $\beta_k >\gamma_k$ for all $k\not\in J$ and $k\neq i$.
Hence, $v(f_i)$ is an absolute point of $\Gamma$. It is immediate
that $v(f_i)$ is irreducible, because $f_i$ is irreducible.

Now, if $\gamma =v(g)\not\in\mathbb{N}^r$ is an irreducible absolute
point of $\Gamma$. Setting $K=I\setminus I_{\gamma}$ we conclude that
$g=h\cdot\prod_{k\in K}f_k^{\alpha_k}$ and $\gamma
=v(g)=v(h)+\sum_{k\in K}\alpha_k\cdot v(f_k)$. As $\gamma$ is
irreducible, we must have $v(h)=\underline{0}$, $K=\{k_0\}$ and
$\alpha_{k_0}=1$, that is, $\gamma =v(f_{k_0})$. \cqd

For the curve $Q = \langle XY \rangle $ the only irreducible
absolute points of $\Gamma$ are $\gamma_1 = v(y) = (1, \infty)$ e
$\gamma_2 = v(x) = (\infty, 1)$, that is, $G = \{x, y\}$ is a
minimal Standard Basis for $\mathcal{O}$ and its semiring is $\Gamma
= \{(0,0)\} \cup \{(1,1) + (\alpha_1, \alpha_2);\ (\alpha_1,
\alpha_2)\in \overline{\mathbb{N}}^2 \}.$ Notice that any element
$(\beta_1, \beta_2) \in \Gamma\setminus\{\underline{\infty}\}$ is
obtained as
$$(\beta_1, \beta_2) = \text{min}\{\beta_1(1, \infty),
\beta_2(\infty, 1)\} = \gamma_1^{\beta_1} \oplus
\gamma_2^{\beta_2}.$$

%-novo
\begin{example} Let us consider the analytic plane curve given by $Q=\langle Y^4-2X^3Y^2-4X^5Y+X^6-X^7\rangle\cup\langle Y^2-X^3 \rangle$ (see Example 3 in \cite{AGMT}).

According Delgado (see \cite{delgado}) the semigroup $S$ of $Q$ is determined by the maximal points
$\{(0,0),(4,2),(6,3),(8,4),(10,5),(12,6),(14,7),(16,8),(18,9),(20,10),(24,12),(22,11),(28,14)\}$
and the semigroup of each branch, that is, $\langle 4,6,13\rangle$ and $\langle 2,3\rangle$.

The set $G$ described in \cite{AGMT} is $\{(4,2),(6,3),(13,15),(26,15),(29,13)\}$, so $S$ is determined by this set and the conductor $(29,15)$ of $S$.

By the other hand, using the above results, the semiring $\Gamma$ of $Q$ is minimally generated by $\{(4,2),(6,3),(13,\infty),(\infty,13)\}$, consequently $S=\Gamma\cap \mathbb{N}^2$.
\end{example}

\vspace{1cm}

\begin{tabular}{lcl}
Carvalho, E. & & Hernandes, M. E. \\
{\it emilio.carvalho$@$gmail.com} & & {\it mehernandes$@$uem.br} \\
 & DMA-UEM &  \\
 & Av. Colombo 5790 & \\
 & Maring\'{a}-PR 87020-900 & \\
 & Brazil &
\end{tabular}

\end{document}